\documentclass{article}

\usepackage{graphicx} 
\usepackage[table]{xcolor}
\usepackage{natbib}
\usepackage{tikz}

\usetikzlibrary{graphs}
\usepackage{url}
\usepackage{hyperref}
\usepackage{amsmath,amsthm,amssymb}
\usepackage{a4wide}
\usepackage{dsfont}
\usepackage{todonotes}
\usepackage{url}
\usepackage{enumitem}
\usepackage{placeins}
\usepackage{float}
\usepackage{authblk}
\usepackage{geometry}
\usepackage{csquotes}
\newtheorem{Theorem}{Theorem}

\newtheorem{cor}{Corollary}

\newtheorem{bsp}{Example}

\newtheorem{Remark}{Remark}

\newtheorem{Def}{Definition}

\newcommand{\I}{\mathbb{I}}
\newcommand{\R}{\mathbb{R}}

\newcommand{\Prob}{\mathbb{P}}

\newcommand{\ignore}[1]{}
\def\multiset#1#2{\ensuremath{\left(\kern-.3em\left(\genfrac{}{}{0pt}{}{#1}{#2}\right)\kern-.3em\right)}}

\def\keywords{\vspace{.5em}
{\textit{Keywords}:\,\relax%
}}
\def\endkeywords{\par}

\title{Extending Characterizations of Multivariate Laws via Distance Distributions}
\author[1]{A. Betken}
\author[1]{A. Marjanovic}
\author[1]{K. Proksch}
\affil[1]{EEMCS, University of Twente, Drienerlolaan 5, 7522 NB Enschede, Netherlands}
\date{}

\begin{document}

\maketitle
\begin{abstract}
We extend a theorem of Maa, Pearl, and Bartoszyński, which links equality of interpoint distance distributions to equality of underlying multivariate distributions, beyond the restrictive class of homogeneous, translation-invariant distance functions. Our approach replaces geometric assumptions on the distance with analytic conditions: volume-regularity of distance-induced balls, Lebesgue differentiability with respect to the distance, and bounded centered oscillations of densities. Under these conditions, equality of interpoint distance distributions continues to imply equality of the generating laws. The result persists under monotone continuous transformations of homogeneous, translation-invariant distances, recovering the original statement, and it extends to compact Riemannian manifolds equipped with the geodesic metric. We further develop a quantitative version of the theorem, i.e. inequalities that connect discrepancies of interpoint distance distributions to the $L^2$-distance between densities, and obtain explicit rates under Ahlfors $\alpha$-regularity of the distance function and  $\beta$-Hölder continuity of densities, capturing dependence on dimensionality. Several representative examples illustrate the applicability of the generalization to domain-specific distances used in modern statistics. The examples include non-homogeneous non-translation invariant distances such as Canberra, entropic distances and the Bray–Curtis dissimilarity. 
\end{abstract}

\section{Introduction} \label{sec::intro}
\noindent
Interpoint distance-based approaches to two-sample testing have been widely used due to their applicability in high-dimensional settings. A variety of methods rely on comparing interpoint distance distributions between and within samples. These include the energy distance \cite{szekely2004testing}, the maximum mean discrepancy (MMD) based test by Gretton et al. \cite{gretton2006kernel}, \cite{gretton2012kernel}, and graph-based statistics such as the minimum spanning tree test and nearest-neighbor methods \cite{Schilling01091986}, \cite{henze1988multivariate}, \cite{bhattacharya2021asymptoticdistributiondetectionthresholds}. Related distance-based procedures include the test of Rosenbaum \cite{rosenbaum2005exact}, which forms an optimal non-bipartite matching of the pooled sample under interpoint distances and tests homogeneity by counting cross-sample pairs.  Montero-Manso and Vilar \cite{montero2019two} proposed a test based on comparing the distributions of interpoint distances directly by  a Cramér–von Mises–type statistic. Rank-based statistics on the interpoint distance distributions have been developed by Liu et al. \cite{liu2022generalized}.  More recently, Betken, Marjanovic and Proksch \cite{betken2024two} introduced a two-sample averaged Wilcoxon test, which applies rank-based statistics to all pairwise distances. The utility of comparing interpoint distance distributions has been further demonstrated in the field of generative modeling, where Jajeśniak et al. \cite{jajesniak2025deep} developed the Interpoint Inception Distance to evaluate models based on their feature representations. These methods share the  structure of reducing the multivariate problem to a univariate comparison of distance-based quantities. \\
The basis for such approaches is provided by a theorem of Maa, Pearl and Bartoszynski \cite{Maa96}, which establishes that, under suitable regularity conditions, the equality of within-sample and between-sample interpoint distance distributions implies the equality of the underlying multivariate distributions. We compile the original statement with all its assumptions in the following theorem. 
\begin{Theorem}[cf. Theorem 2 in \cite{Maa96}]\label{thm:maa}
Let $X_1, X_2, X_3$  be i.i.d. random vectors with values in $\mathbb{R}^k$ and Lebesgue probability density $f_X$, let $Y_1, Y_2, Y_3$ be i.i.d. random vectors with values in $\mathbb{R}^k$ and  Lebesgue probability density $f_Y$ and let $X_1, X_2,X_3$ and $Y_1, Y_2, Y_3$ be independent. Let $d:\mathbb{R}^k\times \mathbb{R}^k\longrightarrow \mathbb{R}$ be a nonnegative, continuous function with
\begin{enumerate}
    \item[\textbf{(D1)}] $d(x,y)= 0$ if, and only if,  $x=y$,
    \item[\textbf{(D2)}] for all $a \in \R$ and $x,y,b \in \R^k$ $d(ax+b, ay+b) = |a|d(x,y)$.
\end{enumerate} 
\noindent
Moreover, assume that
\begin{enumerate}[start=1,label={(\bfseries R\arabic*):}]
    \item[\textbf{(R1)}]  $\int_{\mathbb{R}^k}f^2_X(x) dx, \int_{\mathbb{R}^k}f^2_Y(y)dy<\infty$,
     \item[\textbf{(R2)}]  the zero vector is a Lebesgue point of the function $u(y)= \int f_Y(x+y) f_X(x) dx$, i.e. it holds that $\frac{1}{\lambda(B_r(0))} \int_{B_r(0)} |u(y) - u(0)| d y \underset{r \rightarrow 0}{\rightarrow} 0$, where $B_r(x)$ denotes the ball in $\mathbb{R}^k$ with radius $r$ around $x$ with respect to  the distance function $d$.
\end{enumerate}
\noindent
Then, it holds that
\begin{equation*}
    f_X = f_Y \textit{\quad if, and only if,  \quad } d(X_1, X_2) \overset{\mathcal{D}}{=} d(Y_1, Y_2) \overset{\mathcal{D}}{=}d(X_3, Y_3).
\end{equation*}
\end{Theorem}
\noindent
This theorem reduces a multivariate two-sample problem to a comparison of three induced univariate interpoint distance distributions and serves as a theoretical foundation for distance-based testing, feature selection procedures etc. The assumptions can be categorized into two types:  assumptions on the distance function, i.e., conditions on the topology and associated volume behavior induced by the distance ($\mathbf{(D1)}$-$\mathbf{(D2)}$) and the regularity assumptions on the densities relative to this distance based structure ($\mathbf{(R1)}$-$\mathbf{(R2)}$) .\\
\noindent
However, many metrics, distance functions and dissimilarities used in practice are neither homogeneous nor translation-invariant, and thus fail to satisfy condition $\mathbf{(D2)}$. Examples include the Canberra distance, whose usefulness has been demonstrated in two-sample testing and clustering methods on genomic data \cite{betken2024two}, \cite{proksch2023personalised}, the Bray-Curtis dissimilarity used in microbiology and ecology \cite{ricotta2017some}, and geodesic distances on Riemannian manifolds \cite{chu2024manifold}. These examples illustrate that many practically relevant distances fall outside of the scope of the original Theorem \ref{thm:maa}, which motivates our intention to generalize the theoretical basis to accommodate non-homogeneous, non-translation-invariant distances. \\

\begin{figure}[H]
    \centering
    \includegraphics[width=0.9\textwidth]{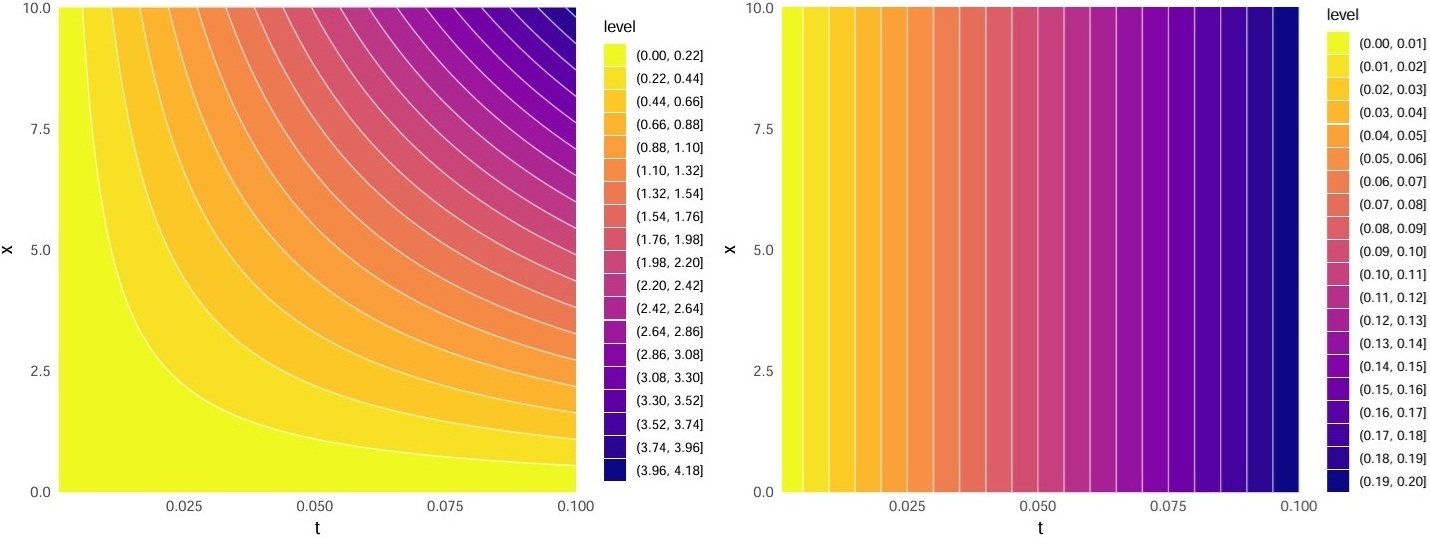}
    \caption{Comparison of the volumes of one-dimensional balls in Canberra and Euclidean distances. The plots represent the volumes $\Phi(x,t)$ as level sets by colorscheme, where $x \in (0,10)$ and $t \in (0, 0.1)$. Left Canberra, right Euclidean balls. For a fixed radius $t$, the volume of an Euclidean ball deos not depend on its position in space, while  in Canberra distance balls become larger the further their center from the origin.  }
    \label{fig:canberra_euclid}
\end{figure}

\begin{figure}[H] 
    \centering
    \includegraphics[width=0.9\textwidth]{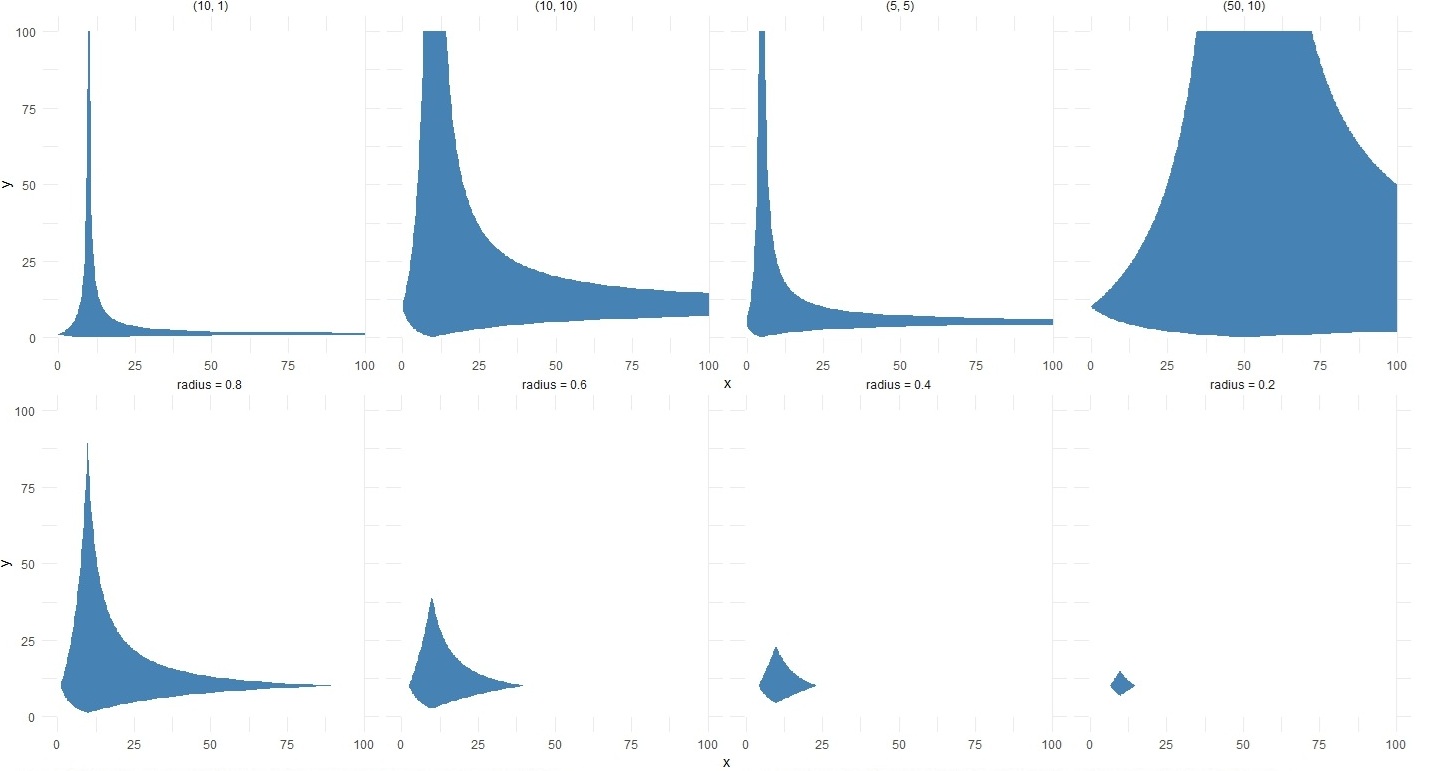}
    \caption{Behavior of two-dimensional balls in Canberra distance. The figures are generated by computing Canberra distances on a dense grid of points and coloring those inside a ball of a given radius. Upper row: unit balls, with centers  $(10, 1)$,   $(10, 10)$,  $(5, 5)$ and  $(50, 10)$. Lower row: balls of radius $0.8$, $0.6$, $0.4$ and $0.2$ centered at $(10, 10)$. }
    \label{fig:canberra_balls}
\end{figure}

\noindent
The proof of the identifiability theorem in \cite{Maa96} leverages homogeneity of the space and polynomial convergence of ball volumes as $t \rightarrow 0$. Assuming translational invariance of $h$, i.e. assumption $\mathbf{(D2)}$ of Theorem \ref{thm:maa} for $a=1$, we have 
    \begin{equation*}
        \mu(B_t^{(h)} (x)) = \mu(B_t^{(h)} (y)) \textnormal{ \quad \quad  } \forall x,y \in \R^k,
    \end{equation*}
    where $B_t^{(h)}(x)$ denotes the ball with center in $x$ and radius $t$ with respect to the distance $h$, while $\mu$ denotes the Lebesgue measure, and therefore $\Phi(x,t) = \Phi(y,t) = \Phi(t)$, where $\Phi(x, t):= \mu(B_t^{(h)} (x)) $. By homogeneity, i.e. assumption $\mathbf{(D2)}$ of Theorem \ref{thm:maa} for $b=0$, it furthermore holds that
        \begin{equation*}
        \Phi(t) =     \mu(B_t^{(h)} (0))=     \mu(|t|B_1^{(h)} (0))=|t|^k  \mu(B_1^{(h)} (0)).
            \end{equation*}
For the Canberra distance on $\R^k$, defined by
\begin{equation}\label{eqn::Canberra_k}
h(x, y) := \sum_{i=1}^k
\begin{cases}
\frac{|x_i - y_i|}{|x_i| + |y_i|}, & \text{if } |x_i| + |y_i| \neq 0 \\
0, & \text{if } x_i = y_i = 0
\end{cases}
\end{equation}
neither of the two conditions, translation invariance and homogeneity, of assumption $\mathbf{(D2)}$ in Theorem \ref{thm:maa} are fulfilled, but the volumes of balls $B_t^{(h)}(x)$ with center in $x$ and radius $t$  do converge to zero as $t$ approaches zero and have comparable convergence rates for different centers in space. To see this, consider the case $k=1$, for which it is $h(x,y) = \frac{|x - y|}{|x| + |y|}$ for $|x| + |y| \neq 0$. In this case, it can be directly calculated that \begin{equation} \label{eqn::Phi_Canberra_1D}
    \Phi(x, t) = \frac{4t \lvert x \rvert}{1-t^2} \textnormal{ \quad \quad for } t \in [ 0,1) \textnormal{ and } x \neq 0.
\end{equation} 
For detailed calculations and expansion to case $k > 1$ see Example \ref{bsp::Canberra}.
Ball volumes in Canberra distance hence behave as rational functions in radius $t$ and are furthermore dependent on the positions of the balls. For visual comparison of one-dimensional Canberra ball volumes to ball volumes of the standard Euclidean distance and  an illustration of the shape of two-dimensional unit balls in Canberra distance see Figures \ref{fig:canberra_euclid} and \ref{fig:canberra_balls}. This motivates the definition of \textit{volume-regular distances} (see Definition \ref{def:volume_reg}), for which class we generalize Theorem \ref{thm:maa} under mild conditions on the densities $f$ and $g$. 


\noindent
In this article, we extend Theorem \ref{thm:maa} beyond homogeneous and translation invariant distance functions: under certain volume and oscillation controls of balls in the considered distance function and assuming Lebesgue differentiability of densities with respect to this distance function, equality of the within- and between-sample distance distributions ensures equality of laws. We complement this finding with computable $L^2$-error bounds expressed via Kolmogorov discrepancies of the three distance distributions. The key idea is to replace the geometric assumptions on the distance, which guarantee the well-behavedness of the balls uniformly in space as the radius shrinks, with more direct analytic conditions on the volumes of the balls.  We provide sufficient criteria under which the equivalence of interpoint distance distributions implies equality of the underlying distributions, and show that these conditions are satisfied for a wide range of distance functions of practical interest. Several corollaries and examples are provided to illustrate the scope of the generalization. The article is organized as follows. Section 2 introduces the standing assumptions (volume-regularity of $h$-balls, $h$-Lebesgue differentiability, bounded centered oscillation), states the generalization, and derives corollaries that recover the original theorem and extend the result to compact manifolds.  Furthermore, in Section 2 we complement the qualitative result with quantitative stability bounds and, under Ahlfors $\alpha$-volume growth and $\beta$-Hölder regularity of the densities, derive dimension-aware rates. In Section 3 we collect illustrative cases, i.e. domain-specific distances, and show that these generate volume-regular balls. Lastly, we outline open questions and future research directions in Section \ref{sec:conclusion}. \\

\section{Generalization}
\noindent
To formulate our generalization of Theorem \ref{thm:maa}, we first specify what we mean by a generalized distance function and the associated notion of balls and their volumes. This isolates the minimal structural ingredients on which our analytic conditions will later be imposed.
\begin{Def}(Generalized distance function, ball and volume) \label{def::gendist}\\
   A function $h:\mathbb{R}^k\times\mathbb{R}^k\to[0,\infty)$ is called  generalized distance function if it is positive definite, i.e. satisfies the identity of indiscernibles:
\begin{align*}
    &h(x,y)= 0 \Leftrightarrow x = y\quad \quad \forall x,y\in\mathbb{R}^k.
\end{align*} 
The ball with center $x$ and radius $t$ with respect to $h$ (or $h$-ball) is the set 
\begin{equation*}
    B_t^{(h)}(x) := \{y \in \R^k|h(x,y) < t\}
\end{equation*}
and the corresponding volume function $\Phi: \R^k \times [0, \infty ) \rightarrow [0, \infty )$ with respect to  the Lebesgue measure $\mu$ on $\mathbb{R}^k$ is defined by $\Phi(x,t) := \mu(B_t^{(h)}(x))$.
\end{Def}
\noindent
Symmetry of the distance function $h$ is not necessary for the proof of our main generalization (Theorem \ref{thm::Main}), but it is additionally assumed for the development of the bounds in Theorem \ref{thm::inequality}. Hence, $h$-balls according to Definition \ref{def::gendist} are, strictly speaking, right-balls.\\
To extend Theorem \ref{thm:maa} beyond homogeneous, translation-invariant metrics, we need to control how volumes of balls behave for generalized  distance functions. In particular,
we find that the volume $\Phi$ of  $h$-balls must shrink  uniformly across different centers in space as their  radii tend to zero. This motivates the notion of \textit{volume-regularity}, which provides the analytic regularity of distance-induced balls needed for generalization of Theorem \ref{thm:maa}.
\begin{Def}(Volume-regular distance function)  \\
A generalized distance function is called volume-regular if it satisfies the following conditions. For some $\epsilon > 0$, the volume function $\Phi$ is  well-defined on $\R^k \times [0,\epsilon)$, i.e. small radius $h$-balls are measurable. 
\begin{equation} \label{eqn::Phi_limit}
     \lim_{t \rightarrow 0^{+}} \Phi(x,t) = 0,
\end{equation}
and for some $\epsilon > 0$ and for all $0 < t < \epsilon$,  $\Phi$ is strictly positive almost everywhere, i.e. 
\begin{equation} \label{eqn::Phi_positive}
     0<\Phi(x,t) <\infty. 
\end{equation}
Lastly, there exists at least one  point $y\in\mathbb{R}^k$ with corresponding $\epsilon > 0$ and $\delta: \R^k \times \R^k \rightarrow (0,\infty)$, such that  for almost every $x\in\mathbb{R}^k$  
\begin{equation} \label{eqn::Phi_uniform}
    c \delta(x,y) \leq  \delta_t (x,y) \leq C \delta(x,y)
\end{equation}
where $ \delta_t (x,y) := \frac{\Phi(x,t)}{ \Phi (y,t)}$ for some constants $c, C > 0$ and for all $0< t < \epsilon$.

\label{def:volume_reg}
\end{Def}
\noindent
We note that the volume function  $\Phi$ is well-defined for (on $\R^k \times \R^k$) measurable $h$. We further note that if condition \eqref{eqn::Phi_uniform} is fulfilled for at least one $y\in\mathbb{R}^k$, then almost any $\xi \in \R^k$ fulfills  \eqref{eqn::Phi_uniform}.  In particular, volume-regularity describes a space in which volumes of $h$-balls centered at different points are stably comparable as radii converge to zero and relative volume $\Phi(x,t)/\Phi(y,t)$ settles for small radii to a finite, positive $\delta (x,y)$ (and is dominated by $\delta (x,y)$) in almost every $x$, while allowing for a controlled divergence of $\Phi$ and $\delta_t$ for a fixed $t$ as $x \rightarrow \infty$. 
Classical Ahlfors-regular geometries, for example, satisfy this condition, and hence imply volume-regularity. We recall the definition of Ahlfors $\alpha$-regularity.
\begin{Def}(Ahlfors $\alpha$-regularity) \label{def::Ahlfors}\\
A metric measure space $(X, h, \mu)$ is called (locally) Ahlfors $\alpha$-regular, if there exists an $\alpha > 0$, $c, C > 0$ and $\epsilon > 0 $, such that 
\begin{equation*}
    c t^\alpha \leq \mu (B^{(h)}_t(x)) \leq C t^\alpha
\end{equation*}
for any $0 < t < \epsilon$ and $x \in X$.
\end{Def}
\noindent
We further relax assumptions on $h$ to generalize Ahlfors $\alpha$-regularity for distance functions as in Definition \ref{def::gendist} and call the distance function $h$ itself Ahlfors $\alpha$-regular (with respect to $\R^k$ and Lebesgue measure).  If $h$ is Ahlfors $\alpha$-regular, then $\frac{c}{ C} \leq \delta_t(x,y) \leq \frac{C}{c}$, i.e. condition \eqref{eqn::Phi_uniform} is fulfilled. The $l_p$-induced distances are, for example, both volume- and Ahlfors-regular. The Canberra distance function on the other hand is volume-regular while not being Ahlfors-regular (see Example \ref{bsp::Canberra}). 
We introduce Definition \ref{def:volume_reg}  as a generalization of Ahlfors regularity, accommodating for  distance functions for which the volume $\Phi$ may not have a polynomial dependence on the radius and may depend on the center of the ball.  \\
\noindent
One further analytic condition needed for the proof of Theorem \ref{thm:maa} is the Lebesgue differentiability relative to the homogeneous translation-invariant distance.  For a generalization of the theorem  we introduce the following concept of Lebesgue differentiability of a locally integrable function $f$ with respect to a function $h$:
\begin{Def} (Lebesgue differentiability with respect to $h$) \\
We call $f \in L^1_{loc} (\R^k)$ Lebesgue differentiable (with respect to $h$) in $x$ if
\begin{equation*}
    \lim_{t \rightarrow 0^{+}} \frac{\int_{B_t^{(h)}(x)} |f(y) - f(x)| d y}{ \Phi(x,t) } = 0.
\end{equation*}
\end{Def}
\noindent
It is important to note that Lebesgue differentiability with respect to a general distance function $h$ does not automatically follow from the assumption $f \in L^1_{\text{loc}}(\mathbb{R}^k)$. In the classical Lebesgue differentiation theorem, differentiability almost everywhere is guaranteed due to the regularity of Euclidean balls and the doubling property of the Lebesgue measure. Lebesgue differentiability with respect to $h$ requires for the underlying (metric) measure space $(\mathbb{R}^k, h, \mu)$ to satisfy conditions such as the doubling property, or more generally, the Vitali covering property. \\
A generalization of the classical differentiation theorem to such settings can be found e.g. in Theorem 1.8 of \cite{heinonen2001lectures}, which ensures differentiability of a locally integrable $f$ on a metric measure space with doubling measure, i.e. such that it holds 
\begin{equation*}
    \mu (B_{2t}^{(h)} (x)) \leq C \mu(B_t^{(h)} (x))
\end{equation*}
for some $C \geq 1$ and any $x$ and $t$, and all balls have finite and positive measure.
\begin{Theorem}(Lebesgue differentiation theorem) \label{thm::Lebesgue_Diff}\\
Let $f \in L^1_{loc}(X)$ where $(X,h,\mu)$ is a doubling metric measure space. Then $f$ is Lebesgue differentiable with respect to $h$, i.e. 
\begin{equation*}
    \lim_{t \rightarrow 0^{+}} \frac{1}{\mu(B_t^{(h)} (x))} \int_{B_t^{(h)} (x)} |f(x) - f(y)| d \mu(y) = 0
\end{equation*}
for almost every $x \in X$.
\end{Theorem}
\noindent
More generally, for a generalized distance $h$ and measure $\mu$, it is necessary that the family of $h$-balls $\{ B_t(x)^{(h)}\}$ forms a Vitali differentiation basis for locally doubling measure $\mu$ for an analogue statement to Theorem \ref{thm::Lebesgue_Diff} to hold  (cf. \cite{folland1999real}). 
\noindent 
Another assumption needed for generalization of Theorem \ref{thm:maa} guarantees  control of the magnitude of  fluctuations of densities within shrinking $h$-ball. This assumption requires the following definition of (uniformly bounded) centered oscillation.
\begin{Def} (Centered oscillation) \\
We call the functional
\begin{equation*}
    A_f^{(h)}(x, \epsilon) := \sup_{0< t < \epsilon} \frac{1}{\Phi(x,t)} \int_{B_t^{(h)}(x)} |f(y) - f(x)| d y
\end{equation*}
centered oscillation of $f$ in $x$ on the scale $\epsilon$ (with respect to $h$). We call $  A_f^{(h)}(x, \epsilon)$ uniformly bounded (on the scale $\epsilon$) if there exist a $C > 0$ such that
\begin{equation*}
     A_f^{(h)}(x, \epsilon) < C \textnormal{ \quad for almost every } x.
\end{equation*}

\end{Def}
\noindent
Uniformly bounded centered oscillation is a mild regularity condition. In fact, many classes of densities, including bounded continuous functions, Lipschitz and Hölder functions,  satisfy it.\\
\noindent
The following theorem extends Theorem \ref{thm:maa} to volume-regular distance functions allowing for
Lebesgue differentiability of the data-generating densities with uniformly bounded, centered oscillations. More precisely, the theorem establishes that equality of interpoint distance distributions implies equality of the underlying distributions under these  assumptions. 
\begin{Theorem} \label{thm::Main}
Let $h:\mathbb{R}^k\times\mathbb{R}^k\to[0,\infty)$ be a volume-regular generalized distance function. Let $X_1$, $X_2$ and $X_3$ be i.i.d. random variables  with density $f$, and $Y_1$, $Y_2$ and $Y_3$ i.i.d. random variables with density $g$, where $f,g$ 
are Lebesgue differentiable with uniformly bounded centered oscillations with respect to $h$ and $X_i$ and $Y_j$ are pairwise independent for all $i, j$.  Let further $f\cdot\delta(\cdot,\xi), g\cdot\delta(\cdot,\xi), f^2\cdot\delta(\cdot,\xi), g^2\cdot\delta(\cdot,\xi)\in L^1(\mathbb{R}^k)$ for a point $\xi\in\mathbb{R}^k$. Then it holds that
\begin{equation} \label{eqn::Maa_Thm}
    f = g \textnormal{ \quad \quad  } \Leftrightarrow \textnormal{ \quad \quad  }  h(X_1, X_2) \overset{\mathcal{D}}{=} h(Y_1, Y_2) \overset{\mathcal{D}}{=} h(X_3, Y_3).
\end{equation}
\end{Theorem}
\vspace{1em}
\noindent
\textbf{Proof. }
\noindent
The forward direction of the equivalence (\ref{eqn::Maa_Thm}) is straightforward, hence it remains to prove the converse. \\
Let $h(X_1, X_2) \overset{\mathcal{D}}{=} h(Y_1, Y_2) \overset{\mathcal{D}}{=} h(X_3, Y_3)$. In particular, then it is 
\begin{equation} \label{eqn::Ps}
     \Prob (h(X_1,X_2) < t) =  \Prob (h(Y_1,Y_2) < t) =  \Prob (h(X_3,Y_3) < t)
\end{equation}
for any $t \in \R$.
Since by assumption (\ref{eqn::Phi_uniform}) for $t$ small enough it uniformly holds
\begin{equation*}
    f^j(x)\delta_t(x, \xi) \leq C f^j(x)\delta(x, \xi),
\end{equation*}
the integrability of $f^j \cdot \delta_t(\cdot, \xi)$ for $j = 1, 2$ follows from  $ f^j\cdot\delta(\cdot,\xi) \in L^1(\R^k)$. Furthermore, since
\begin{equation*}
    \int_{\R^k} f^j(x) \Phi(x,t) dx = \Phi(\xi, t) \int_{\R^k} f^j(x) \delta_t(x, \xi) dx
\end{equation*}
the integrability of $f^j \cdot \Phi(\cdot, t)$ follows. The analog holds for $g^j \cdot \Phi(\cdot, t)$, $g^j \cdot \delta_t(\cdot, \xi)$. \\ 
\noindent
It holds
\begin{align*}
    \Prob (h(X_1,X_2) < t) &= \int_{\R^k} \int_{B_t^{(h)} (x)} f(y) f(x) dy dx \\
    &=\int_{\R^k} f(x) \int_{B_t^{(h)}(x)} (f(y) - f(x)) d y dx  + \int_{\R^k} f^2(x) \Phi(x,t) dx.
\end{align*}
Therefore,
\begin{equation} \label{eqn::pre_limit}
    \frac{\Prob (h(X_1,X_2) < t)}{\Phi(\xi, t)} = \int_{\R^k} f(x) \frac{\int_{B_t^{(h)}(x)} (f(y) - f(x))dy}{\Phi(\xi, t)}  dx  + \int_{\R^k} f^2(x) \delta_t(x,\xi) dx.
\end{equation}
Invoking the Lebesgue differentiability  of $f$ with respect to $h$ it holds for almost every $x$ that
\begin{equation}\label{eq:pwc}
    \lim_{t\to0^+}f(x) \delta_t (x, \xi) \frac{\int_{B_t^{(h)}(x)} (f(y) - f(x)) dy}{\Phi(x, t) } = 0.
\end{equation}
Since the centered oscillation of $f$ is uniformly bounded and (\ref{eqn::Phi_uniform}) holds, we find
\begin{align*}
    \left|f(x) \delta_t (x, \xi) \frac{\int_{B_t^{(h)}(x)} (f(y) - f(x)) dy}{\Phi(x, t) }\right|\leq A_f^{(h)}(x, \epsilon) f(x)\delta_t(x,\xi) \leq C f(x)\delta(x,\xi)
\end{align*}
for a.e. $x$ with some constant $C>0$ and some $\epsilon > 0$, for all $0 < t < \epsilon$.
By point-wise convergence \eqref{eq:pwc} and integrability of the majorant $f\cdot \delta(\cdot,\xi)$, an application of Lebesgue's dominated convergence theorem yields
\begin{equation} \label{eqn::dom_convergence}
 \lim_{t \rightarrow 0^{+}} \int_{\R^k} f(x) \frac{\int_{B_t^{(h)}(x)} (f(y) - f(x))dy}{\Phi(\xi, t)}  dx   = 0.
\end{equation}
Hence from (\ref{eqn::pre_limit}) it follows
\begin{equation*}
   \lim_{t \rightarrow 0^{+}} \frac{\Prob (h(X_1,X_2) < t)}{\Phi(\xi, t)} =  \lim_{t \rightarrow 0^{+}}  \int_{\R^k} f^2(x) \delta_t(x,\xi) dx.
\end{equation*}
Analogously, we obtain
\begin{equation*}
   \lim_{t \rightarrow 0^{+}} \frac{\Prob (h(Y_1,Y_2) < t)}{\Phi(\xi, t)} =  \lim_{t \rightarrow 0^{+}}  \int_{\R^k} g^2(x) \delta_t(x,\xi) dx
\end{equation*}
and
\begin{equation*}
   \lim_{t \rightarrow 0^{+}} \frac{\Prob (h(X_3,Y_3) < t)}{\Phi(\xi, t)} = \lim_{t \rightarrow 0^{+}}   \int_{\R^k} f(x)g(x) \delta_t(x,\xi) dx
\end{equation*}
where $fg \cdot \delta(\cdot, \xi) \in L^1(\R^k)$ follows with Cauchy-Schwarz inequality and implies integrability of $fg \cdot \delta_t(\cdot, \xi)$.\\

\noindent
From (\ref{eqn::Ps}) it hence follows
\begin{align*} 
    0 &=  \lim_{t \rightarrow 0^{+}} \frac{ \Prob (h(X_1,X_2) < t) +  \Prob (h(Y_1,Y_2) < t) - 2 \Prob (h(X_3,Y_3) < t)  }{ \Phi(\xi, t)} \\
    & = \lim_{t \rightarrow 0^{+}} \int_{\R^k} (f(x) - g(x))^2 \delta_t(x,\xi) d x.
\end{align*}
Since it is a.s. $\delta_t(x,\xi) \geq c \delta(x,\xi) > 0$ by assumption (\ref{eqn::Phi_uniform}),  we conclude f = g almost everywhere. \qed 
    
\begin{Remark}
As in the original theorem (cf. proof of Theorem 2 in \cite{Maa96}), full independence of the $X$– and $Y$–samples is not required; it suffices that $X_{1}$ and $X_{2}$, $Y_{1}$ and $Y_{2}$, and $X_{3}$ and $Y_{3}$ are pairwise independent.
In the original setting, the assumptions of homogeneity and translational invariance of $h$ implied that the $h$-ball volume $\Phi(x,t)$ was independent of the center point $x$. This allowed the proof to be stated in terms of a convolution of $f$ and $g$ around zero. In Theorem \ref{thm::Main}, these geometric assumptions are replaced by the assumption of Lebesgue differentiability with respect to $h$.
Moreover, the role of homogeneity in controlling the asymptotic behavior of ball volumes (as $t\rightarrow 0$) is taken over by conditions  (\ref{eqn::Phi_limit}) and (\ref{eqn::Phi_uniform}), which ensure that the volume ratio $\Phi(x,t)/\Phi(\xi,t)$ does not diverge and that balls shrink uniformly  as $t \rightarrow 0$.
Finally, since $\Phi(x,t)$ may now depend on $x$ and need not be constant for fixed $t$, additional integrability assumptions are required. These control divergence of $\Phi(x,t)$ as $x \rightarrow \infty$, a phenomenon that indeed occurs for certain distances such as the Canberra and entropic distance (see Examples \ref{bsp::Canberra} and \ref{bsp::entropic}). 
\end{Remark}
\begin{Remark}
The ratio $\delta_t$ does not necessarily converge (to some $\delta$) as $t \rightarrow 0$. To see this, consider the following example of a volume-regular distance: Let $h: \R^k \times \R^k \rightarrow [0,\infty)$ with $h(x,y):= g_x(|x-y|)$,  where  for some  $t_0 > 0$, $g_x: (0, r_x(t_0)) \rightarrow (0, t_0)$ is the inverse function of 
\begin{equation*}
    r_x(t):= t\left( 1 + \epsilon A(x) \sin\left(\log\left(\frac{1}{t}\right)\right)\right)
\end{equation*}
  where $\epsilon A(x)$ with $A(x) < 1$ is a small amplitude dependent on $x$ and $\epsilon$ small enough. In particular, $g_x$ is well-defined, since $r_x'(t) = 1 + \epsilon A(x) (\sin(\log(\frac{1}{t})) -\cos(\log(\frac{1}{t}))) > 0$ for $0 < t < t_0$ and hence $r_x(t)$ is strictly monotone increasing with $r_x (t) \underset{ t \rightarrow 0^+}{ \rightarrow} 0$. Since for $0 < t < t_0$
  \begin{equation*}
      B^{(h)}_t(x) = \{y \in \R^k | h(x,y) <  t \} = \{ y \in \R^k | \quad |x-y| < r_x(t)  \} =  B^{(l_1)}_{r_x(t)}(x)
  \end{equation*}
  we therefore have
  \begin{equation*}
      \Phi(x,t) = \mu( B^{(l_1)}_{r_x(t)}(x)) = \omega_{1,k} r_x^k(t)
  \end{equation*}
  where $\mu$ denotes Lebesgue measure on $\R^k$ and $\omega_{1,k}$  the volume of $k$-dimensional unit ball in $l_1$. Since it is $(1-\epsilon)t  <  r_x(t) <  (1+\epsilon)t $ for all $x \in \R^k$ and $0 < t < t_0$, it follows
  \begin{equation*}
      \omega_{1,k} (1-\epsilon)^k t^k  < \Phi(x,t)  < \omega_{1,k} (1+\epsilon)^k t^k
  \end{equation*}
  and $h$ hence volume-regular, while 
  \begin{equation*}
      \delta_t(x,y) = \left( \frac{1 + \epsilon A(x) \sin(\log(\frac{1}{t}))}{1 + \epsilon A(y) \sin(\log(\frac{1}{t}))} \right)^k.
  \end{equation*}
  Having $A(x) \neq A(y)$, $ \delta_t$ does not converge for $t \rightarrow 0$: it is, however, bounded for any $x, y$ and $0 < t < t_0$.  \\
\end{Remark}
\noindent 
\cite{Maa96} noted that their result should extend to monotone transformations of the Euclidean metric; see Remark 2 in \cite{Maa96}. Corollary \ref{thm::corollary_1} verifies this statement for monotone continuous transformations of homogeneous, translation-invariant distances.  In particular, for $\gamma = \textnormal{id}$ we recover the classical theorem.

\begin{cor} \label{thm::corollary_1}
Let $X_1$, $X_2$ and $X_3$ be i.i.d. with density $f$, $Y_1$, $Y_2$ and $Y_3$ i.i.d. with density $g$, where $f,g \in L^2(\R^k)$ and $X_i$ and $Y_j$ independent for all $i, j$. Let $d$ be a homogeneous and translation invariant distance function on $\R^k$, let $\gamma: [0, \infty) \rightarrow [0, \infty)$  be continuous and strictly increasing with $\gamma(0) = 0$ and $h := \gamma \circ d $. If $f$ and $g$ have uniformly bounded centered oscillation (with respect to $h$), then

\begin{equation*}
    f = g \textnormal{ \quad \quad  } \Leftrightarrow \textnormal{ \quad \quad  }  h(X_1, X_2) \overset{\mathcal{D}}{=} h(Y_1, Y_2) \overset{\mathcal{D}}{=} h(X_3, Y_3).
\end{equation*}

\end{cor}
\noindent
\textbf{Proof. }
Since strictly monotone increasing functions  are invertible, it holds that  
\begin{equation*}
    B_t^{(h)}(x) = \{y \mid h(x,y) < t \} = \{y \mid d(x,y) < \gamma^{-1} (t)\}=  B_{ \gamma^{-1} (t)}^{(d)}(x).
\end{equation*}
It follows that
\begin{equation*}
    \mu (B_t^{(h)}(x)) = \mu (B_{\gamma^{-1}(t)}^{(d)}(x)).
\end{equation*}
For a homogeneous and translation-invariant distance function $d$ it further holds that 
 \begin{equation*}
     \mu (B_{r}^{(d)}(x)) = r^k \mu (B_{1}^{(d)}(0))
 \end{equation*}
and we conclude
\begin{equation*}
    \Phi (x,t) = \left(\gamma^{-1}(t) \right)^k \mu (B_{1}^{(d)}(0))
\end{equation*}
It follows with $\gamma^{-1}(t) \rightarrow 0$ for $t \rightarrow 0^{+}$ that $\Phi$ fulfills the conditions (\ref{eqn::Phi_limit}), (\ref{eqn::Phi_positive}) and (\ref{eqn::Phi_uniform}).\\
 The integrability conditions of Theorem \ref{thm::Main} follow from   $f,g \in L^2(\R^k)$, since the volume in $h$ is translational-invariant, i.e. $\Phi$ does not depend on $x$. The Lebesgue differentiation limit with respect to $h$ is the same as along $d$ and since $(\R^k, d, \mu)$ is doubling, the Lebesgue differentiability follows from $f, g \in L^1(\R^k)$. With uniformly bounded centered oscillations of $f$ and $g$, the corollary follows from Theorem \ref{thm::Main}.  \qed \\
 
\noindent
A natural question is whether our result extends to curved spaces, where translation invariance and homogeneity fail. Compact Riemannian manifolds provide a prototypical example: their geodesic balls have controlled volume growth and the Riemannian measure is doubling. Corollary \ref{cor::manifold} illustrates that our generalization continues to hold in this broader geometric setting. 
\begin{cor} \label{cor::manifold}
    Let $\mathcal{M}$ be a smooth, compact $k$-dimensional Riemannian manifold equipped with Riemannian measure $\mu$ and geodesic distance $d_{\mathcal{M}}$. Let on $\mathcal{M}$ be given i.i.d. r.vs. $X_1$, $X_2$ and $X_3$ with density $f$ and i.i.d. $Y_1$, $Y_2$ and $Y_3$ with density $g$, where  $f,g \in L^2(\mathcal{M})$ . Then it  holds that
    \begin{equation*}
    f = g \textnormal{ \quad \quad  } \Leftrightarrow \textnormal{ \quad \quad  }  d_{\mathcal{M}}(X_1, X_2) \overset{\mathcal{D}}{=} d_{\mathcal{M}}(Y_1, Y_2) \overset{\mathcal{D}}{=} d_{\mathcal{M}}(X_3, Y_3).
\end{equation*}
\end{cor}
\noindent
\textbf{Proof. }
The geodesic ball volume on $\mathcal{M}$ satisfies 
\begin{equation*}
    \Phi(x,t)= \omega_{2,k}\, t^k\!\left[
1 - \frac{\mathrm{S}(x)}{6(k+2)}\,t^2
+ \frac{-3\|R\|^2(x) + 8\|\mathrm{Ric}\|^2(x) + 5\,\mathrm{S}(x)^2 - 12\,\Delta \mathrm{S}(x)}
{360\,(k+2)(k+4)}\,t^4
+ O(t^6)
\right]
\end{equation*}
where $\omega_{2,k}$ denotes the volume of the $k$-dimensional Euclidean unit ball, $S$ denotes scalar curvature, $R$ Riemannian curvature tensor and  $\mathrm{Ric}$ Ricci tensor (cf. Theorem 3.1 in \cite{gray1974volume}). Conditions (\ref{eqn::Phi_limit}), (\ref{eqn::Phi_positive}) and (\ref{eqn::Phi_uniform}) follow. Since $\mathcal{M}$ is compact, it is also
\begin{equation} \label{eqn::Phi_bound_corollary}
    c t^k < \Phi(x,t) < Ct^k  \textnormal{ \quad \quad in a.e. } x
\end{equation}
for some $c, C >0$ and integrability conditions of Theorem \ref{thm::Main} follow from $f, g \in L^2(\mathcal{M})$. Furthermore, $\mathcal{M}$ is with (\ref{eqn::Phi_bound_corollary}) a metric measurable space of globally doubling measure and by Lebesgue differentiation theorem, $f$ and $g$ are Lebesgue differentiable with respect to $d_{\mathcal{M}}$.\qed\\

\noindent
Theorem \ref{thm::Main} shows that equality of the interpoint distance distributions implies equality of the underlying distributions, but it remains qualitative. In practice, one is also interested in quantifying how close the densities are when the interpoint distributions are only approximately equal. To this end, we aim to establish inequalities that connect the $L^2$–distance between $f$ and $g$ with Kolmogorov-type distances between probability functions $F_{XX}$, $F_{XY}$, $F_{YY}$, where $F_{XX}(t) :=  \Prob (h(X_1, X_2) < t)$, $F_{XY}(t) :=  \Prob (h(X_3, Y_1) < t)$ and $F_{YY}(t) :=  \Prob (h(Y_2, Y_3) < t)$. Such bounds allow us to measure the deviation of $f$ and $g$ in terms of observable discrepancies of interpoint distances, and provide a quantitative version of Theorem \ref{thm::Main}.

\begin{Theorem}  \label{thm::inequality}
Let $h$ be a symmetric generalized distance function, the conditions of Theorem \ref{thm::Main} be fulfilled and furthermore let for $\delta$ hold 
\begin{equation} \label{eqn::delta_bonds}
    0 < \delta_*  < \delta(x,\xi) < \delta^* < \infty  \textnormal{\quad \quad}  \textnormal{\quad for a.e. } x \in \R^k
\end{equation}
for some $\xi \in \R^k$ and  $\delta_*, \delta^*$. Let $\Delta_K(t)$ denote Kolmogorov discrepancy, i.e.
\begin{equation*}
    \Delta_K(t) := \sup_{0 < u \leq t} \left| F_{XX} (u) - F_{XY} (u)\right| + \sup_{0 < u \leq t} \left| F_{YY} (u) - F_{XY} (u)\right|.
\end{equation*}
For some $\epsilon > 0 $ and $0 < t < \epsilon$ it holds then
\begin{equation} \label{eqn::ineq_L2}
    ||f-g||^2_{L^2} \leq \frac{1}{c \delta_*} \left[ \frac{\Delta_K(t)}{\Phi(\xi , t)} + r(\xi, t) \right]
\end{equation}
and
\begin{equation} \label{eqn::ineq_delta}
    \Delta_K (t) \leq C \delta^* \Phi(\xi, t)  (||f||_{L^2} + ||g||_{L^2}) ||f-g||_{L^2},
\end{equation}
where $c$ and $C$ are as in Definition \ref{def:volume_reg},  with $r(\xi, t) \rightarrow 0$ as $t \rightarrow 0$.
\end{Theorem}
\noindent
\textbf{Proof.}
Following (\ref{eqn::pre_limit}), analogous integrals for $\Prob(h(X_3, Y_1) < t)$ and $\Prob(h(Y_2, Y_3) < t)$ are
\begin{align*}
    \frac{F_{XY}(t)}{\Phi(\xi, t)} &= \frac{1}{\Phi(\xi, t)} \int_{\R^k} f(x) \int_{B_t^{(h)}(x)}(g(y) - g(x) ) d y d x + \int_{\R^k} f(x) g(x) \delta_t(x, \xi) dx \\
    &= \frac{1}{\Phi(\xi, t)} \int_{\R^k} g(x) \int_{B_t^{(h)}(x)}(f(y) - f(x) ) d y d x + \int_{\R^k} g(x) f(x) \delta_t(x, \xi) dx
\end{align*}
due to the symmetry  of $h$ and the first integral, and
\begin{equation*}
    \frac{F_{YY}(t)}{\Phi(\xi, t)} =  \frac{1}{\Phi(\xi, t)} \int_{\R^k} g(x) \int_{B_t^{(h)}(x)}(g(y) - g(x) ) d y d x + \int_{\R^k}  g^2(x) \delta_t(x, \xi) dx
\end{equation*}
and therefore
\begin{align*}
    \int_{\R^k} (f(x) - g(x))^2 \delta_t (x, \xi) d x &= \frac{1}{\Phi (\xi, t)} \left( F_{XX}(t) - F_{XY}(t)+F_{YY}(t) - F_{XY}(t) \right) \\
    &- \frac{1}{\Phi (\xi, t)}  \int_{\R^k} (f - g)(x) \int_{B_t^{(h)} (x)} \left((f - g)(y) - (f -g)(x)  \right) d y d x  \\
    & \leq \frac{\Delta_K (t)}{\Phi (\xi, t)} + r(\xi, t)
\end{align*}
where 
\begin{equation*}
    r(\xi, t):= \frac{1}{\Phi (\xi, t)}  \int_{\R^k} |(f -g)(x)| \int_{B_t^{(h)} (x)} \left|(f - g)(y) - (f - g)(x) \right| d y d x .
\end{equation*}
In particular, since uniformly bounded centered oscillations of $f$ and $g$ imply uniformly bounded centered oscillation of $(f-g)$ (i.e. with $A_{|f-g|} \leq A_f + A_g$), by dominated convergence analogously to the proof of Theorem \ref{thm::Main} it follows that $r(\xi, t) \rightarrow 0$ as $t \rightarrow 0$. 
Finally, by assumption  there exist $\epsilon > 0 $ s.t. for $0 < t < \epsilon$ and a.e. $x$ it is $c \delta_* \leq c \delta(x, \xi) \leq \delta_t(x, \xi)$ and hence
\begin{equation*}
  c  \delta_* ||f-g||^2_{L^2} \leq  \int_{\R^k} (f(x) - g(x))^2 \delta_t (x, \xi) d x,
\end{equation*}
and the inequality (\ref{eqn::ineq_L2}) follows. \\

\noindent
On the other hand, as we have
\begin{equation*}
    F_{XX}(u) - F_{XY}(u) = \int_{\R^k} f(x) \int_{B_u^{(h)} (x)}  (f(y) - g(y)) d y d x
\end{equation*}
by applying Cauchy-Schwarz's inequality, it follows that
\begin{equation} \label{eqn::CS1}
    |F_{XX}(u) - F_{XY}(u)| \leq \left( \int_{\R^k}   f^2(x)d x \right)^{\frac{1}{2}}  \left( \int_{\R^k}  \left( \int_{B_u^{(h)} (x)}  (f(y) - g(y)) d y \right)^2 d x \right)^{\frac{1}{2}}
\end{equation}
and
\begin{align*}
     \int_{B_u^{(h)} (x)}  (f(y) - g(y)) d y & = \int_{\R^k}  (f(y) - g(y)) \I^2_{B_u^{(h)}(x)}(y) d y \\
     & \leq \left(\int_{\R^k} \I_{B_u^{(h)}(x)}(y) d y  \right)^{\frac{1}{2}}  \left(\int_{\R^k}(f(y) - g(y))^2 \I_{B_u^{(h)}(x)}(y) d y  \right)^{\frac{1}{2}} \\
     &\leq \Phi(x, u)^{\frac{1}{2}} \left(  \int_{B_u^{(h)} (x)}  (f(y) - g(y))^2 d y \right)^{\frac{1}{2}}.
\end{align*}
\noindent
Since it is $\Phi(x,u) = \delta_t(x, \xi) \Phi(\xi, u ) $ and $\delta_t(x, \xi) \leq C\delta^*$ for small $t$ by assumption, it holds

\begin{align*}
 \int_{\R^k}  \left( \int_{B_u^{(h)} (x)}  (f(y) - g(y)) d y \right)^2 d x & \leq \int_{\R^k}  \Phi(x,u)  \int_{B_u^{(h)} (x)}  (f(y) - g(y))^2 d y d x \\
 & \leq C \delta^* \Phi(\xi, u)  \int_{\R^k} \int_{\R^k} (f(y) - g(y))^2 \I_{B_u^{(h)} (x)} (y) dy dx.
\end{align*}
Since $h$ is symmetric and therefore $\I_{B_u^{(h)} (x)} (y) = \I_{B_u^{(h)} (y)} (x)$, by Fubini it is 
\begin{align*}
    \int_{\R^k} \int_{\R^k} (f(y) - g(y))^2 \I_{B_u^{(h)} (x)} (y) dy dx & =    \int_{\R^k} (f(y) - g(y))^2 \int_{\R^k} \I_{B_u^{(h)} (y)} (x) d x d y \\
    & = \int_{\R^k} (f(y) - g(y))^2 \Phi(y, u) dy \\
    & \leq C \delta^* \Phi(\xi, u) ||f-g||^2_{L^2}
\end{align*}
\noindent
Combining the inequalities we find
\begin{equation*}
    |F_{XX}(u) - F_{XY}(u)| \leq C\delta^* \Phi(\xi, u) ||f||_{L^2}||f-g||_{L^2}.
\end{equation*}
Similarly, by considering the integral
\begin{equation*}
     F_{YY}(u) - F_{XY}(u) = \int_{\R^k} g(x) \int_{B_u^{(h)} (x)} (g(y) - f(y)) d y d x
\end{equation*}
in this case analogously we conclude
\begin{equation*}
      |F_{YY}(u) - F_{XY}(u)| \leq C\delta^* \Phi(\xi, u) ||g||_{L^2}||f-g||_{L^2}.
\end{equation*}
Taking the supremum over $0<u\leq t$ and using  monotonicity of $\Phi(\xi,u)$ in $u$ we obtain 
\begin{equation*}
    \Delta_K (t) = \sup_{0 < u \leq t} \left| F_{XX} (u) - F_{XY} (u)\right| + \sup_{0 < u \leq t} \left| F_{YY} (u) - F_{XY} (u)\right| \leq C\delta^* \Phi(\xi, t) ( ||f||_{L^2} + ||g||_{L^2})||f-g||_{L^2}
\end{equation*}
proving (\ref{eqn::ineq_delta}). \qed 
\begin{Remark}
The derived inequalities hold locally for some $\epsilon > 0$ and $0 < t < \epsilon$, stemming from the scale of volume-regularity and boundedness of centered oscillations. Furthermore, applying Theorem \ref{thm::inequality} in order to  bound the Kolmogorov distance, is not informative since $\Phi(\xi, t)$ generally diverges as $ t \rightarrow \infty$. Nevertheless, a global reverse bound can be obtained by combining (\ref{eqn::ineq_delta}) on the interval $u \leq t$ with a separate tail estimate for $u > t$, e.g. via moment, Markov-type bounds, which we will not pursue here. 
\end{Remark}
\noindent
The bounds on $\delta(x, \xi)$ simplify comparability of $||f-g||_{L^2}$ with the integral over $(f-g)^2 \delta_t$. For standard distances like $l_p$-norms or geodesic distances on compact manifolds, condition (\ref{eqn::delta_bonds}) is automatically true, while this does not hold for the  Canberra distance, the entropic distance or the Bray-Curtis dissimilarity and local versions of inequalities (\ref{eqn::ineq_L2}) and (\ref{eqn::ineq_delta}) have to be developed. 
\noindent
Additionally, the inequalities of Theorem \ref{thm::inequality} are general but involve the quantities $\Phi(x,t)$ and $r(\xi,t)$, whose behavior is not explicit in full generality. By imposing additional structural assumptions on both the distance function $h$ and the densities $f$ and $g$, these terms can be controlled more directly. In particular, Ahlfors $\alpha$–regularity of $h$ provides precise asymptotics for the volume $\Phi(x,t)$, while Hölder continuity of $f$ and $g$ ensures that the remainder term $r(\xi, t)$ can be bounded uniformly. Under these combined conditions, the inequalities of Theorem \ref{thm::inequality} simplify to the more concrete statement given in the following corollary. In particular, under the assumption of  Ahlfors $\alpha$–regularity of $h$, the upper bound in Theorem \ref{thm::inequality} is proportional to $\Delta_K^{\tfrac{\beta}{\alpha+\beta}}$, where $\Delta_K = \Delta_K(\infty)$. Accordingly, as $\alpha$ increases, this upper bound deteriorates. Since $\alpha$ typically describes the ambient dimension of the data, this upper bound reflects that quantitative recovery of $f$ and $g$ from interpoint distances becomes weaker in higher dimensions. \\

\begin{cor} \label{cor::dim_corollary}
Let the conditions of Theorem \ref{thm::inequality} be fulfilled and furthermore let $h$ be Ahlfors $\alpha$-regular. Let $f$ and $g$ be $\beta$-Hölder continuous with respect to $h$, i.e. there exist $ C_f, C_g > 0 $ s.t.  
 \begin{equation*}
     |f(x) - f(y)| \leq C_f (h(x,y))^\beta  \textnormal{, \quad} |g(x) - g(y)| \leq C_g (h(x,y))^\beta  \textnormal{ \quad } \forall x, y \in \R^k.
 \end{equation*}
Then it holds
\begin{equation} \label{eqn::ineq_kolmogorov}
    ||f-g||^2_{L^2} \leq C \Delta_K^{\frac{\beta}{\alpha + \beta}}
\end{equation}
where $\Delta_K = \Delta_K (\infty)$.
\end{cor}
\noindent
\textbf{Proof.}
We have for some $\epsilon > 0 $ and $c_*, c^* > 0$ that $c_* t^\alpha \leq \Phi(x,t) \leq c^* t^\alpha$ for all $x\in \R^k$ and $0 < t  < \epsilon$. Applying Theorem \ref{thm::inequality} we have
\begin{align*}
     ||f-g||^2_{L^2} &\leq \frac{1}{\delta_*} \frac{\Delta_K(t)}{ c_* t^\alpha} + \frac{1}{\delta_*} \frac{1}{c_* t^\alpha}  \int_{\R^k}f(x) \int_{B_t^{(h)} (x)} \left|f(y) - g(y) - (f(x) - g(x)) \right| d y d x \\ 
     &+ \frac{1}{\delta_*} \frac{1}{c_* t^\alpha}\int_{\R^k}g(x) \int_{B_t^{(h)} (x)} \left|f(y) - g(y) - (f(x) - g(x)) \right| d y d x  \\
     &\leq \frac{1}{\delta_*} \frac{\Delta_K}{ c_* t^\alpha} +  \frac{1}{\delta_*} \frac{1}{c_* t^\alpha} \int_{\R^k}(f(x)+g(x)) \Phi(x,t) ( \sup_{y \in B_t^{(h)}(x)} |f(y) - f(x)| + \sup_{y \in B_t^{(h)}(x)} |g(y) - g(x)|) dx \\
     & \leq \frac{1}{\delta_*} \frac{\Delta_K}{ c_* t^\alpha} +   \frac{1}{\delta_*} \frac{1}{c_* t^\alpha} c^* t^{\alpha} (\sup_{x \in \R^k} \sup_{y \in B_t^{(h)}(x)} |f(y) - f(x)| + \sup_{x \in \R^k} \sup_{y \in B_t^{(h)}(x)} |g(y) - g(x)|)  \int_{\R^k}(f(x)+g(x)) dx .
\end{align*}
By Hölder continuity there exist $c$ such that for any $x$, $y$ it is $|f(x) - f(y)| \leq c h(x,y)^\beta$ and $|g(x) - g(y)| \leq c h(x,y)^\beta$. With $h(x,y) < t $ we conclude
\begin{align*}
    ||f-g||^2_{L^2} \leq \frac{1}{\delta_*} \frac{\Delta_K}{ c_* t^\alpha} + 2\frac{1}{\delta_*}\frac{c^*}{c_*} c t^\beta.
\end{align*}
Minimization of $\phi(t)= \frac{1}{\delta_*} \frac{\Delta_K}{ c_* t^\alpha} + 2\frac{1}{\delta_*}\frac{c^*}{c_*} c t^\beta$ over $t$ yields
\begin{equation*}
    t_0 = ( \frac{\alpha \Delta_K}{2 \beta c^* c})^{\frac{1}{\alpha + \beta}}
\end{equation*}
and it follows
\begin{equation*}
     ||f-g||^2_{L^2} \leq C \Delta_K^{\frac{\beta}{\alpha + \beta}}
\end{equation*}
where $C$ depends on $\alpha$, $\beta$, $\delta_*$, $c_*$, $c^*$ and $c$.
\qed

\begin{Remark}
\noindent
    As an example, consider $\R^k$ equipped with the $l_p$-distance, in which case it holds
    \begin{equation*}
        \Phi(x,t) = \Phi(t) = \omega_{p,k} t^k,
    \end{equation*}
    where $\omega_{p,k}$ is the unit ball volume. Let $f$ and $g$ be Lipschitz continuous. Then, Corollary \ref{cor::dim_corollary} implies that 
    \begin{equation*}
         ||f-g||^2_{L^2} \leq C \Delta_K^{\frac{1}{k + 1}}.
    \end{equation*}
\end{Remark}

\section{Examples}
The general conditions of Theorem \ref{thm::Main}, i.e. volume-regularity, can be checked directly for many distances used in applications. In this section we work through several representative cases. Starting with norm-induced distances, we then consider Canberra and entropic distances, both of which do not fulfill the homogeneity and translation invariance assumption of Maa’s theorem, before concluding with Bray–Curtis dissimilarity, a weighted version of the Canberra distance. These examples illustrate that the considered framework covers a broad range of practically relevant settings. 
\begin{bsp}($\| \cdot \|_p^p$ induced distance)
Let $h(x,y) =\| x-y \|_p^p $. As the composition of $(\cdot)^p$ and the $l_p$-distance, $h$ satisfies the conditions of Corollary \ref{thm::corollary_1}. If the conditions on $f$ and $g$ are fulfilled, (\ref{eqn::Maa_Thm}) follows. In particular, $h$ becomes a component-wise additive quantity, in contrast to $l_p$-distance.
\end{bsp}
\noindent
Next we consider the Canberra distance, a distance function often applied in clustering and used in high-dimensional applications. It is widely applied in ecology, \enquote{omics}-counts and spectrometry, in situations when the data features are sparse and it is necessary to emphasize relative differences.  While we exactly calculated the volume of a Canberra-ball in one dimension (see Section \ref{sec::intro}), we
provide upper and lower bounds in higher dimensions
proving  volume-regularity. 
\begin{bsp} (Canberra distance) \label{bsp::Canberra}\\ 

\noindent
Let us consider the Canberra distance of $x=(x_1, \ldots, x_k)^{\top}$ and $y=(y_1, \ldots, y_k)^{\top}$ in $\R^k$ defined by
\begin{equation}\label{eqn::Canberra_k}
h(x, y) := \sum_{i=1}^k
\begin{cases}
\frac{|x_i - y_i|}{|x_i| + |y_i|}, & \text{if } |x_i| + |y_i| \neq 0, \\
0, & \text{if } x_i = y_i = 0.
\end{cases}
\end{equation}
\noindent
\textbf{Case $k = 1$:} \\
In this case, $h(x,y) = \frac{\lvert x - y \rvert}{\lvert x\rvert +  \lvert y\rvert}$. Consider for $x > 0$ and $t \in (0,1)$ the cases $y \geq x$ and $y  < x$ separately. 
For $y \geq x$ and   $h(x,y)  < t$, we have
\begin{equation*}
    \frac{y - x}{x + y}  < t, \ \text{or, equivalently,} \  \   y  < \frac{x(t + 1)}{1 - t}.
\end{equation*}
For  $y  < x$ and   $h(x,y)  < t$, we have
\begin{equation*}
     \frac{x-y}{x + y}  < t, \ \text{or, equivalently,} \  \   y > \frac{x(1-t)}{1+t}.
\end{equation*}
It follows that
\begin{equation} \label{eqn::Canberra_1D_cube}
    B_t(x)= \left(x\frac{1-t}{1+t} , x\frac{1+t}{1-t} \right)
\end{equation}
and therefore $\mu(B_t(x)) = \frac{4tx}{1-t^2}$. Due to symmetry, for $x < 0$ we have $\mu(B_t(x)) = -\frac{4tx}{1-t^2}$ and we conclude
\begin{equation} \label{eqn::Phi_1D}
    \Phi(x, t) = \frac{4t \lvert x \rvert}{1-t^2} \textnormal{ \quad \quad for } t \in (0,1) \textnormal{ and } x \neq 0.
\end{equation}
Moreover, it holds that  
\begin{equation*}
    B_t(0) = 
\begin{cases}
\{ 0\} & t < 1 , \\
\R & t \geq 1,
\end{cases}
\end{equation*}
and therefore 
\begin{equation*}
     \Phi(0,t) = 
\begin{cases}
0 & t < 1 , \\
\infty & t \geq 1.
\end{cases} 
\end{equation*}

\noindent
\textbf{Case $k > 1$:} \\
We firstly provide an upper and lower bounds on the volume of $B_t^{(h)}(x)$ for the case $x_i \neq 0$ for all $i \in \{1,...,k\}$. For the upper bound, we consider the rectangle in $\R^k$ with sides corresponding to  \eqref{eqn::Canberra_1D_cube}. Since $y_i \notin \left(x_i\frac{1-t}{1+t} , x_i\frac{1+t}{1-t} \right) $ implies $h(x,y) \geq t$, we have
\begin{equation*}
    B_t(x) \subseteq \prod\limits_{i=1}^k \left(x_i\frac{1-t}{1+t} , x_i\frac{1+t}{1-t} \right).
\end{equation*}
It follows that 
\begin{equation} \label{eqn::Canberra_upper}
    \Phi(x,t) \leq \left(\frac{4t}{1-t^2} \right)^k \prod\limits_{i=1}^k \lvert x_i \rvert \textnormal{ \quad \quad for } t \in (0,1) \textnormal{ and } x_i \neq 0.
\end{equation}
For the lower bound consider the rectangle with sides 
\begin{equation*}
   I_i = \left(x_i\frac{1-t/k}{1+t/k} , x_i\frac{1+t/k}{1-t/k} \right).
\end{equation*}
Since $y_i \in I_i$ for all $i \in \{1,...,k\}$ implies $y \in B_t(x)$, we conclude that
\begin{equation*}
     \prod\limits_{i=1}^k \left(x_i\frac{1-t/k}{1+t/k} , x_i\frac{1+t/k}{1-t/k} \right) \subseteq B_t(x)
\end{equation*}
and therefore
\begin{equation}\label{eqn::Canberra_lower}
    \left(\frac{4t/k}{1-(t/k)^2} \right)^k \prod\limits_{i=1}^k \lvert x_i \rvert \leq \Phi(x,t) \textnormal{ \quad \quad for } t \in (0,1) \textnormal{ and } x_i \neq 0.
\end{equation}
If for any $i \in \{1,...,k\}$ $x_i = 0$, then in that dimension the set of permissible $y_i$ is empty, and therefore $\Phi(0,t) = 0$ for $t \in (0,1)$. The set $\{x=(x_1, \ldots, x_k)^{\top} \mid x_i = 0 \ \text{for at least one $i \in \{1, \ldots, k\}$}\}$ has, however, Lebesgue measure $0$ in $\R^k$.
Accordingly,  for a.e. $x$
\begin{equation*}
    \Phi(x,t) \sim C(x) \left( \frac{t}{1-t^2} \right)^k,
\end{equation*}
where $0<C(x)\sim \prod_{i=1}^k \lvert x_i \rvert$. More precisely, with \eqref{eqn::Canberra_upper} and \eqref{eqn::Canberra_lower} it follows that $\delta_t(x,y) \leq C \delta(x,y)$ with $\delta(x,y) = \prod_{i= 1}^k \left|\frac{x_i}{y_i}\right|$ for some $C>0$. We conclude that the conditions \eqref{eqn::Phi_limit}, \eqref{eqn::Phi_positive} and \eqref{eqn::Phi_uniform} in the definition of volume-regular functions, i.e. Definition \ref{def:volume_reg},  are fulfilled.  Theorem \ref{thm::Main} then applies if its conditions on the densities $f$ and $g$ are fulfilled.\\

\noindent
As a concrete example, consider the density function $f$ of the standard normal distribution, i.e. $f(x) = \frac{1}{\sqrt{2 \pi}} e^{-\lVert x \rVert^2 /2}$, which is smooth on $\R^k$. Then, the Lebesgue differentiability condition is fulfilled since
\begin{align*}
      \frac{1}{\Phi(x,t)} \int_{B_t(x) } \lvert f(x) - f(y)\rvert d y &\leq \frac{1}{\Phi(x,t)} \sup_{y \in B_t(x)}\lvert f(x) - f(y)\rvert \Phi(x,t) \\
      &= \sup_{y \in B_t(x)}\lvert f(x) - f(y)\rvert \rightarrow 0
\end{align*}
for almost every $x \in \R^k$ and $t \rightarrow 0^{+}$. According to \eqref{eqn::Canberra_upper}, it further holds that
\begin{equation*}
    f^2(x)\Phi(x,t) \leq \frac{1}{2\pi}\left(\frac{4t}{1-t^2} \right)^k \prod\limits_{i=1}^k \lvert x_i \rvert  e^{-\lVert x \rVert^2 }
\end{equation*}
with the right-hand side of the above inequality being  integrable on $\R^k$. We conclude that $f \delta(\cdot,y), f^2 \delta(\cdot,y) \in L^1(\R^k)$ for any $t \in (0,1)$ and the integrability condition is fulfilled. In other words: Theorem \ref{thm::Main} holds for  normally distributed data and the Canberra distance. 
\end{bsp}
\noindent
Note that for the Canberra distance $\Phi$ and $\delta$  are not globally bounded and diverge as $x \rightarrow \infty$, hence do not fulfill condition \eqref{eqn::delta_bonds} of Theorem \ref{thm::inequality}. This can, however, be amended by locally applying Theorem \ref{thm::inequality} on a bounded set (i.e. such that $\delta$ is bounded) carrying the majority of the mass of $f$ and $g$, while  separately developing bounds for the tails. A detailed discussion of the latter, however, lies beyond the scope of this work.\\
\noindent
We further treat the entropic distance, which  involves logarithmic terms, impeding exact volume calculation even in one dimension. Instead, we provide upper and lower bounds on the entropic volume in one dimension and then proceed to generalize these to higher dimensions. 
\begin{bsp} (Entropic distance on $\R^k_{>0}$) \label{bsp::entropic} 
\noindent
Consider $h: \R^k_{>0}\times \R^k_{>0} \rightarrow \left[ 0, \infty \right)$ defined by
\begin{equation*}
    h(x,y):= \sum_{i=1}^k \left| x_i \log \left( \frac{x_i}{y_i} \right) -x_i + y_i \right|.
\end{equation*}
As in the previous example, we firstly consider the case $k=1$ and subsequently generalize bounds for $\Phi(x,t)$ from one to more dimensions. \\

\noindent
\textbf{Case $k=1$:} 
In this case, $h(x,y) =\left| x \log \left( \frac{x}{y} \right) -x + y \right| $. Let $\eta_x (y) := x \log \left( \frac{x}{y} \right) -x + y  $ s.t. $h(x,y)= |\eta_x (y)|$. 
Since $x>0$, for a ball of radius $t < \frac{x}{2}$ around $x$, let $y \in (x-t, x+t)$ and therefore 
\begin{equation} \label{eqn::I}
    y-x \in \left(-\frac{x}{2}, \frac{x}{2}\right), \ \text{or, equivalently, \ }     y \in \left(\frac{x}{2}, \frac{3x}{2}\right).
\end{equation}
A second-order Taylor expansion  of  $\eta_x (y)$ around $x$ yields
\begin{equation*}
    \eta_x (y) = \frac{1}{2x} (y-x)^2  \frac{1}{3} \frac{x}{\zeta^3} (y-x)^3
\end{equation*}
for some $\zeta \in  \left(\frac{x}{2}, \frac{3x}{2}\right)$. Furthermore, it holds that
\begin{equation} \label{eqn::II}
    |\eta_x (y)| < t \quad \Leftrightarrow \quad \frac{1}{2}(y-x)^2 \left|\frac{1}{x} - \frac{2}{3} \frac{x}{\zeta^3}(y-x) \right|< t.
\end{equation}
With the triangle inequality and \eqref{eqn::II} it follows that
\begin{align*}
    \left|\frac{1}{x} - \frac{2}{3} (y-x) \right| \leq \frac{1}{x} + \frac{2}{3} \frac{x}{\zeta^3 } |y-x|
    & < \frac{1}{x} + \frac{2}{3} \frac{x}{\left( \frac{x}{2} \right)^3 } \frac{x}{2} 
     = \frac{11}{3} \frac{1}{x}
\end{align*}
and
\begin{equation*}
     \left|\frac{1}{x} - \frac{2}{3} \frac{x}{ \zeta^3} (y-x) \right| > \frac{1}{x} - \frac{1}{3} \frac{1}{x} = \frac{2}{3}\frac{1}{x}
\end{equation*}
since $\zeta$ is lower-bounded by $x$ for $y-x > 0 $ and $y-x$ is upper-bounded by $\frac{x}{2}$ according to \eqref{eqn::I}.  
 
\noindent
In particular, for any $y$ satisfying the first inequality, i.e. such that
\begin{equation*}
    y \in (x - \sqrt{\frac{11}{6}xt}, x + \sqrt{\frac{11}{6}xt})
\end{equation*}
it is also $ |\eta_x (y)| < t $, while any $y \in B_t^{(h)}(x)$ satisfies the second inequality and we conclude
\begin{equation*} \label{eqn::k=1}
    \left(x - \sqrt{\frac{11}{6}xt}, x + \sqrt{\frac{11}{6}xt}\right) \subset B_t^{(h)}(x) \subset   \left(x - \sqrt{3xt}, x + \sqrt{3xt}\right).
\end{equation*}
\noindent
For any two points $x$ and $y$ and balls around them with radii $t$ we have
\begin{equation*}
      \lim_{t \rightarrow 0^{+}} \Phi( x,t)\geq \lim_{t \rightarrow 0^{+}} 2\sqrt{\frac{6}{11}tx} = 0
\end{equation*}
and 
 \begin{equation}
    c\frac{\sqrt{x}}{\sqrt{y}}< \frac{    \Phi(x,t) }{    \Phi(y,t) }< C\frac{\sqrt{x}}{\sqrt{y}},
 \end{equation}
such that conditions \eqref{eqn::Phi_limit} - \eqref{eqn::Phi_uniform} are fulfilled. \\

\noindent
\textbf{Case $k>1$:}\\
Similarly to Example \ref{bsp::Canberra}, we construct rectangles in and containing $B_t^{(h)}(x)$ using the intervals derived for case $k=1$. Following the previous derivation of inscribed and circumscribed intervals of the ball $B_t^{(h)}(x)$ in one dimension, we derive

\begin{equation*}
     \Pi_{i=1}^k \left(x_i -\sqrt{\frac{11}{6}x_i \frac{t}{k}} , x_i +\sqrt{\frac{11}{6}x_i \frac{t}{k}}\right) \subset B_t^{(h)}(x) \subset \Pi_{i=1}^k \left(x_i -\sqrt{3x_i t}, x_i +\sqrt{3x_i t})\right) 
\end{equation*}
and hence a lower and an upper bound on $\Phi(x,t)$ are 
\begin{equation*}
  \left(2 \sqrt{\frac{11}{6k}} \right)^k t^{k/2} \Pi_{i=1}^k x_i^{1/2} < \Phi(x,t) < (2 \sqrt{3})^k t^{k/2} \Pi_{i=1}^k x_i^{1/2}.
\end{equation*}
We conclude that
\begin{equation*}
    \Phi(x,t) = C(x) t^{k/2} 
\end{equation*}
where $C(x) > 0$ is a function with behavior similar to $\Pi_{i=1}^k x_i^{1/2}$. Conditions (\ref{eqn::Phi_limit}), (\ref{eqn::Phi_positive}) and (\ref{eqn::Phi_uniform}) hold with $\delta (x,y) = \Pi_{i = 1}^k \left(\frac{x_i}{y_i} \right)^{\frac{1}{2}}$. \\
If for the density function $f$ it holds  $ f \in L^2(\mathbb{R}_{>0}^k) $ and the weighted function $ f^2(x) \left( \prod_{i=1}^k x_i \right)^{1/2} \in L^1(\mathbb{R}_{>0}^k) $, given a decay of $f(x)$ quicker than the growth of $ \left( \prod_{i=1}^k x_i \right)^{1/4}$, then the integrability conditions for such $f$ are fulfilled. \\
Lastly, since on a bounded set (such is any $B_t^{(h)}(x)$) $C(x)$ is bounded, there exist $0< C_1 < C_2$ with $$C_1 t^{k/2} < \Phi(x,t) < C_2 t^{k/2}$$ and since $\Pi_{i=1}^k I_{i, t} \subset B_t^{(h)}(x)$ for any $t$ where $ I_{i, t} := \left(x_i -\sqrt{\frac{11}{6}x_i \frac{t}{k}} , x_i +\sqrt{\frac{11}{6}x_i \frac{t}{k}}\right) $, such $f \in L^1_{loc}(\R^k)$ is Lebesgue differentiable in a.e. $x$ by a generalized Lebesgue differentiation theorem (cf. Theorem 3.21 in \cite{folland1999real}). 
\end{bsp}
\noindent
As for Example \ref{bsp::Canberra}, we note that the bounds of Theorem \ref{thm::inequality} are not directly applicable due to divergence of $\delta$ at infinity. Finally, we consider the Bray–Curtis dissimilarity, used in ecology and microbiome studies. Rather than coordinate-wise normalization as Canberra, Bray-Curtis is globaly normalized, making it less sensitive to rare features. Its close relation to Canberra distance makes the argumentation indirect and we do not derive the estimates on volume as in previous examples. 
\begin{bsp} (Bray-Curtis dissimilarity on $\R^k_{>0}$) \\
The Bray-Curtis dissimilarity on $\R^k_{>0}$ is given by
\begin{equation*}
    h_{BC} (x,y) = \frac{\sum_{i=1}^k |x_i - y_i|}{\sum_{i=1}^k (x_i + y_i)}.
\end{equation*}
Ricotta-Podani (\cite{ricotta2017some}) rewrite $h_{BC}$ as
\begin{equation*}
     h_{BC} (x,y) = \sum_{i=1}^k \omega _i (x,y)\frac{ |x_i - y_i|}{ x_i + y_i}
\end{equation*}
where
\begin{equation*}
    \omega _i (x,y) = \frac{x_i + y_i}{\sum_{i=1}^k (x_i + y_i) }
\end{equation*}
and $\sum_{i=1}^k \omega_i = 1$. Bray-Curtis dissimilarity can thus be seen as weighted Canberra distance where a weight $\omega_i$ is continuous, strictly positive and uniformly bounded between $0$ and $1$. Consequently, the volume of a Bray-Curtis ball satisfies the same asymptotics (as $t \rightarrow 0$) as in the Canberra case and we conclude that the conditions of Theorem \ref{thm::Main} hold.\\
\end{bsp}

\section{Conclusion} \label{sec:conclusion}
We showed that interpoint-distance distribution characterize data distributions under  mild conditions on the distance function and the data-generating distributions. The framework recovers the classical setting of the Maa–Pearl–Bartoszyński Theorem (see Corollary \ref{thm::corollary_1}), extends to compact manifolds (see Corollary \ref{cor::manifold}), and yields quantitative $L^2$-bounds in terms of Kolmogorov discrepancies with dimension-aware rates under Ahlfors-type volume growth and Hölder regularity (Theorem \ref{thm::inequality} and Corollary \ref{cor::dim_corollary}). These rates make explicit the trade-off between the ambient dimension $k$ and the Hölder-regularity index $\beta$, providing interpretable control in high dimensions. The quantitative bounds are presently local and do not apply to Canberra, Bray–Curtis and entropic distances, since their volume ratios $\delta_t(x,y)$ are not bounded in $x \in \R^k$. The obstruction is the blow-up of coordinate-normalized ratios near axes or at the origin, which breaks uniform comparability required by our global arguments. Still, it is possible to extend the results of Theorem \ref{thm::inequality} and Corollary \ref{cor::dim_corollary} locally for these distances as well. A practical route is to decompose the space into a truncated region where ratios are controlled and a tail region handled by distributional mass bounds. 
A global bound should follow by combining this localized estimate based on Corollary \ref{cor::dim_corollary} with additional tail estimates that quantify the mass of $f$ and $g$ outside of the truncated region. \\
\noindent
Beyond this limitation, the results cover many distances used in practice and clarify when multivariate two-sample problems may be reduced to one-dimensional comparisons. We have, in particular, developed a theoretical foundation for testing based on the Canberra distance, as proposed in \cite{betken2024two}. This provides the practical procedure with identifiability guarantees and, on truncated domains, the same dimension-aware rates as for Ahlfors regular spaces.

\bibliography{bibliography}
\bibliographystyle{alpha}

\end{document}